\documentclass[10pt]{amsart}
\usepackage{amsmath,amssymb,enumerate}
\usepackage{epsf,epsfig,amsfonts,graphicx,color}
 \usepackage[applemac]{inputenc}
\usepackage{url}
\usepackage{a4wide,epsfig}
 \numberwithin{equation}{section}
\newtheorem{theorem}{Theorem}

\newtheorem{remark}[theorem]{Remark}
\newcommand{\weg}[1]{}

\newtheorem{proposition}[theorem]{Proposition}
\newtheorem{cor}[theorem]{Corollary}

\title[Generic metric has no nontrivial Killing tensors]{The geodesic flow of a generic metric does not admit
nontrivial integrals polynomial in momenta}
\author{Boris Kruglikov and Vladimir S. Matveev}
 \address{B.\,S. Kruglikov: Department of Mathematics and Statistics,
University of Troms\o, 9037 Troms\o, Norway.
Email: {\tt boris.kruglikov@uit.no}}
 \address{V.\,S. Matveev: Institute of Mathematics, Friedrich-Schiller-Universit\"at, 07737 Jena, Germany.
Email: {\tt vladimir.matveev@uni-jena.de}}
\begin{document}

 \begin{abstract}
Any smooth geodesic flow is locally integrable with smooth integrals. We show that generically this fails if
we require, in addition, that the integrals are polynomial (or, more generally, analytic) in momenta. Consequently
we obtain that a generic real-analytic metric does not admit, even locally, a real-analytic integral.
 \end{abstract}

\maketitle

\vspace{-0.5cm}

\section{Introduction: definitions,  results and motivation.}
Let $g=g_{ij}(x)dx^idx^j$ be a pseudo-Riemannian metric of arbitrary signature (including the Riemannian one)
on a connected manifold $M$ of dimension $n\ge 2$.
Consider the geodesic flow of $g$, which is the Hamiltonian system on the cotangent bundle with
$H= \tfrac12g^{ij}(x)p_ip_j$. Here and below  $x=(x^1,\dots,x^n)$ are local coordinates
on $M$ and $p=(p_1,\dots,p_n)$ are the corresponding coordinates on the cotangent spaces.

A function $F:T^*M\to  \mathbb{R}$ is called an \emph{integral} if $\{F,H\}\equiv 0$, i.e.
 \begin{equation} \label{bracket}
\sum_{i=1}^n \Bigl(\frac{\partial H}{\partial x^i}\frac{\partial F}{\partial p_i} -
\frac{\partial F}{\partial x^i}\frac{\partial H}{\partial p_i}\Bigr)\equiv  0.
 \end{equation}
A function $F$ is called polynomial in momenta of degree $d$ if it has the form
 $$
F(x,p)=\sum_{i_1,\dots,i_d=1}^n  K^{i_1\dots i_d}(x)\,p_{i_1}\cdots p_{i_d},
 $$
where the coefficients $K^{i_1...i_d}(x)$ may depend on  $x$; without loss of generality we assume that
the coefficients  are invariant  with respect to permutations of the indices.

An equivalent terminology used in differential geometry and relativity is \emph{Killing tensors}:
one immediately sees that components $K^{i_1...i_d}$ behave under the coordinate changes as components
of a $(d,0)$ tensor; it is known that \eqref{bracket} is equivalent to the \emph{Killing equation}
 \begin{equation} \label{killing}
K_{(i_1...i_d, i_{d+1})}\equiv 0.
 \end{equation}
Here we use the metric $g$ to lower indices,  denote by  ``comma'' the covariant differentiation in the Levi-Civita connection of $g$, and  use the brackets $( ... )$ to indicate  the symmetrization with respect to  the indices staying inside.

The geometric condition equivalent to equation \eqref{killing} is that for any geodesic $\gamma$,
whose velocity vectors will be denoted by $\dot\gamma^i$, the function
\begin{equation} \label{condkill}
t\mapsto \sum_{i_1,\dots,i_d=1}^n K_{i_1...i_d}(\gamma(t))\,
\dot\gamma(t)^{i_1} \dots \dot\gamma(t)^{i_d} \ \textrm{does not depend on $t$}.
 \end{equation}

We say that an  integral   is \emph{trivial}, if it is a function of $H$. Clearly, a trivial integral which
is polynomial in momenta of odd degree $d$ is identically zero, and trivial integrals of even degree $d=2q$ are
proportional to $H^q$ with a constant coefficient of proportionality.

Notice that integrability of the geodesic flow of $g$ on $M$ (or simply the existence of an integral of $H=H_g$
on $T^*M$) implies the same for a domain $D\subset M$. Our results are actually nonexistence results, so  can suppose $D$ to be a ball in $\mathbb{R}^n$, since   nonexistence of an integral of a generic  metric on  $D$ implies
nonexistence of an integral of a generic  metric on any manifold  $M$. Our goal is to prove the following statements.

 \begin{theorem} \label{thm:1}
Let $g$ be a $C^k$-smooth metric on an open disc $D\subset\mathbb{R}^n$, where $k\ge 2$. Then, for any $d\in\mathbb{N}$ and $\varepsilon>0$ there exists a metric $\tilde g$ on $D$, which is $\varepsilon$-close to $g$ in the $C^k$-topology,
and $\varepsilon'>0$ such that for any $C^2$-smooth metric $g'$ on $D$, which is $\varepsilon'$-close to $\tilde g$
in the $C^2$-topology, the geodesic flow of the metric $g'$ does not admit a nontrivial integral polynomial in momenta of degree $d$.
 \end{theorem}

In Theorem \ref{thm:1} we do not require any smoothness of the (coefficients of the) integral. They need not to be even
continuous (note that to define non-smooth integrals we cannot use formulae \eqref{bracket} or \eqref{killing};
instead we use \eqref{condkill}). This is because as a by-product of the proof of Theorem \ref{thm:1} we obtain:

 \begin{cor}\label{cor:0}
Let $F$ be an integral of the geodesic flow of a metric $g$. Assume that $F$ is polynomial in momenta of degree $d$.
If the metric is $C^k$-smooth for $k\ge 2$, then the integral is  $C^{k-1}$-smooth.
 \end{cor}

We actually conjecture that if the metric is $C^{k\ge 2}$-smooth, then the integral is also $C^k$-smooth.
This conjecture holds for $n=2, d\leq2$, and it can also be confirmed in some other special cases.

Recall that a subset of a topological space is called \emph{generic}, if it is a countable intersection of open
everywhere dense subsets (this implies that such set is of second category in the Baire sense).
Elements of a generic subset are called \emph{generic} elements. We apply this definition to the topological space of
$C^k$-smooth metrics on $D$ (this space is  a Baire space, so generic elements are everywhere dense) and obtain:

 \begin{cor}\label{cor:1}
For $k\ge 2$, generic (in the $C^k$-topology) $C^k$-smooth metrics admit only trivial integrals  polynomial in momenta.
 \end{cor}

Recall that, as it was known at least to Darboux \cite{Darboux} (see also \cite{Wh}), the existence of an
integral that is real analytic in momenta (dependence of  the coefficients on the position can be only smooth)
implies the  existence of an integral that is polynomial in momenta. Thus, we have:

 \begin{cor}\label{cor:2}
For $k\ge 2$, generic (in the $C^k$-topology) $C^k$-smooth metrics admit only trivial integrals analytic in momenta.
 \end{cor}

A very special case of Corollary \ref{cor:2}  (dimension is 2, metric is real analytic) follows from \cite{Ten}.

Let us recall that \emph{$N$-jet  of a function} $f$ at a point $x$ can be identified with the collection
$\bigl(f(x),\tfrac{\partial f}{\partial x_1}(x),\dots,\tfrac{\partial^N f}{\partial x_n^N}(x)\bigr)$
of the values of $f$ and its partial derivatives up to order $N$ at $x$. Similarly,
\emph{jet of a  metric $g$} can be identified with the collection of jets of the entries $g_{ij}$ of the metric considered as functions.

 \begin{theorem} \label{thm:2}
For  arbitrary $n\ge 2$,  $d\ge 1$,  there exists  $N= N(n,d)\in \mathbb{N}$ with the following property:

For any $C^N$-smooth metric $g$ on $D\subset \mathbb{R}^n$, any $x\in D$ and any $\varepsilon >0$
there exists a $C^N$-smooth metric $\tilde g$ on $D$, which is $\varepsilon$-close to $g$ in the $C^N$-topology,
and $\varepsilon'>0$ such that for any $C^N$-smooth metric $g'$ on $D$, whose $N$-jet at $x$ is $\varepsilon'$-close
to that of $\tilde g$, the geodesic flow of the restriction of $g'$ to any open connected
subset $D'\subset D$ containing $x$ does not admit a nontrivial integral polynomial in momenta of degree $d$.
 \end{theorem}

One can show that $N(n,d)$ can be chosen to be
\begin{equation} \label{displayed}
d+1+\tfrac{(n+d-1)!(n+d)!}{(n-1)!n!d!(d+1)!}.
\end{equation}
This $N(n,d)$ is definitely not the minimal value such that Theorem \ref{thm:2} is valid, which is expected to be much smaller. This can be already seen for $d=1$: to explicitly compute the minimal values of $N(n,1)$ observe
that integrals linear in momenta are the same as Killing vector fields
(symmetries of the metric). Using this, one can show that $N_{\min}(2,1)=4$
and $N_{\min}(n,1)=3$ for $n=3,4,\dots$, while the formula \eqref{displayed} gives
$N(2,1)=5$, $N(3,1)=8$, $N(4,1)=12$, $N(5,1)=17$, $N(6,1)=23$.

It is not easy to find the minimal  value for $N(n,d)$ for $d>1$ (such that Theorem \ref{thm:2} remains valid). From \cite{kruglikov} (see also
 \cite{duna}) it follows  $N_{\min}(2,2)=6$, and this result is pretty involved.
From the answer of Bryant on \cite{mathoverflow} it follows that $N_{\min}(2,d)\ge 2d +2$, and in fact we can show that $N_{\min}(2,d)= 2d +2$.
In our opinion, it is an interesting and challenging problem to find, for other $(n,d)$,  the minimal  or at least much  better values of $N(n,d)$ than that given by  \eqref{displayed}. The first open nontrivial case is $N(3,2)$.

Theorems \ref{thm:1} and \ref{thm:2} will not surprise mathematicians treating the Killing equation
via the geometric theory of partial differential equation. Indeed, system of PDEs \eqref{killing}
on the components  $K_{i_1...i_d}$ is overdetermined; the coefficients of this system
depend on the components of the metric and their partial derivatives. One expects therefore that for generic data
(i.e., for a generic  metric) the compatibility conditions do not allow nontrivial solutions. However even the experts
could be puzzled with the $C^2$-smoothness assumptions in Theorem \ref{thm:1}, since the compatibility conditions require  higher  order of derivatives. From this viewpoint Theorem \ref{thm:2} is more expectable.

This actually was our initial approach to this problem, but it appeared that
it is not that easy  to compute the compatibility conditions for all $d$ and $n$.
According to our knowledge  it was never done
even for $(n,d)=(2,3)$ and $(3,2)$, and our calculations show that they could be quite intractable even for
these small values of $n$ and $d$ (see also \cite{kruglikov,Wolf}).
Instead of calculation of the integrability conditions we  invented a  trick  which solves the problem.

For mathematicians with background in dynamical systems and ergodic theory Theorem \ref{thm:1} and especially
Corollary \ref{cor:2} may look less evident. An important open conjecture   in the field (see e.g.  
\cite[\S 8.1]{burns}) is  that {\it every closed manifold (of dimension $\ge 2$)
admits a metric with ergodic geodesic flow}  --- most experts are sure
that this conjecture is true but there is no proof and we are not aware of any  idea that may lead to  a proof in full generality.  Recall that ergodicity of the geodesic flow is equivalent to nonexistence of an integral (which allowed to be merely a  measurable function).

In  dimensions $n\ge 3$ this problem  is   almost completely open.
For two-dimensional manifolds under some additional assumptions on the integral (e.g. real-analyticity),
the conjectures  follows from the combination of results \cite{paternain,ContrerasPa,knieper}:
a generic metric has positive topological entropy and this obstructs integrability in the class of real analytic integrals.
See also \cite{Bolotin} where on any surface it is  constructed  a metric  such that the geodesic flow  of any $\varepsilon$-close, in the $C^2$-topology,  metric does not admit a nontrivial  real analytic integral.

Note that a version of the latter claim holds true in any dimension due to the combination of results \cite{Contreras,taimanov}:
a generic metric on a closed $n$-dimensional manifold does not admit $n$ functionally independent real-analytic
integrals in involution. Nonexistence of a non-trivial single integral for a generic (even real-analytic)
metric of dimension $n\ge 3$ has not been known.

Corollary \ref{cor:2} resolves the above problem in all dimensions $n\ge2$ assuming real analyticity of the
integral in momenta: \emph{a generic metric $g$ on $M$ has no such integral}. Moreover, the result is local
(our manifold $M$ is not assumed to be closed and can be a disc). This local result is wrong
if we merely assume that the integrals are smooth: locally the geodesic flow of any metric is smoothly integrable.
This result is known in the folklore and we formulate it and give a sketch of the proof at the end of the paper, see Proposition \ref{lastt}.

The explanation of this visual contradiction  is that for a generic metric (which can be even real analytic), the integral  (which can be  real analytic on the unit bundle $S_1T^*D=\{H=\frac12\}$) cannot be real-analytically extended to
the whole cotangent bundle; we shall make this point more clear in the Conclusion.

The proofs  are  organized as follows. In \S\ref{sec:2.2} we prove Theorem \ref{thm:1}. As  mentioned above,
the proof is based on a certain trick. To make it easier for the reader, we first explain the trick for $d=1$
in \S \ref{sec:2.1}. Of course, the case $d=1$ of Theorem \ref{thm:1} is mostly evident from the geometric point of view, since integrals with $d=1$ correspond to Killing vector fields. Our proof contains
all ideas of the general case and is easier to digest than the general proof (a hurried reader could go directly to
\S\ref{sec:2.2}). Corollary \ref{cor:0} will be explained in \S\ref{sec:2.3}. Corollaries  \ref{cor:1}, \ref{cor:2} are trivial. Theorem \ref{thm:2} is implied by Theorem \ref{thm:1} and one (folklore) known result that we formulate with
a sketch of the proof  and apply in \S \ref{sec:2.4}. Finally,   in \S \ref{sec:2.5} we will demonstrate local
integrability of any (Riemannian or pseudo-Riemannian) metric with smooth integrals, and indicate further
applications of methods developed in our paper.

\subsection*{Acknowledgements.}  We thank G. Paternain who encouraged us to write this paper, partially because of
the results of this paper were applied in     his paper \cite{preprint}, joint with H. Zhou, for useful discussions and pointing us the reference \cite{ozeki}.
We also thank the user Joe, who keeps an air of mystery about him, for asking the question \cite{mathoverflow} at Mathoverflow, and R. Bryant for active participation in this  Mathoverflow discussion, from which it is clear that he knew  Theorems \ref{thm:2}  and \ref{folklore} at least in dimension $n=2$.
We are also grateful to  D. Treschev for useful discussions and for attracting our interest to  \cite{Bolotin}.
The paper  was written  during and shortly after  the visit of VM to the University of Troms\o\ supported by the DAAD
project  57068692  (PPP Norwegien); VM thanks the University of Troms\o\  for hospitality.

\section{ Proof of Theorem \ref{thm:1} and  its corollaries.} \label{sec:2}

\subsection{ Proof of Theorem \ref{thm:1} in the case $d=1$. } \label{sec:2.1}

We work locally, in a small geodesically convex disc $D\subset \mathbb{R}^n$  with   $n\ge 2$.
Notice that the property of geodesic convexity is stable under $C^2$-small perturbations of the metric $g$.

On the disc $D$ consider $3n$ points, which we denote by $ A_1,...,A_n, B_1,...,B_n, C_1,...,C_n$. We denote the set $\{A_1,...,A_n\}$ by $A$, the  set $\{B_1,...,B_n\}$ by $B$ and  the  set $\{C_1,...,C_n\}$ by $C$.

Assume that no three of the points $ A_1,...,A_n, B_1,...,B_n, C_1,...,C_n$ lie on one geodesic and, in the case of indefinite signature,
 that no geodesic connecting any two of these points is light-like.
We also assume that the points are generic in the following sense: for any point from any  of this set and for any other
set the initial vectors of geodesics connecting that point with all point of  this other set are linearly independent.
(In the proof for general $d$ the condition that the initial vectors are linearly independent will be replaces by a more complicated algebraic condition)

Denote by $S^1_{A_i}$ (resp.\ $S^1_{B_i}$, $S^1_{C_i}$) the space of linear functionals on $T_{A_i}D$
(resp.\ on $T_{B_i}D$, $T_{C_i}D$). It is an $n$-dimensional vector space. We will denote the elements of $S^1_{A_i}$ by $\alpha_i$ and  the elements of $S^1_{B_i}$ by $\beta_i$.

 \begin{figure}[ht!]
{\includegraphics[width=0.7\textwidth]{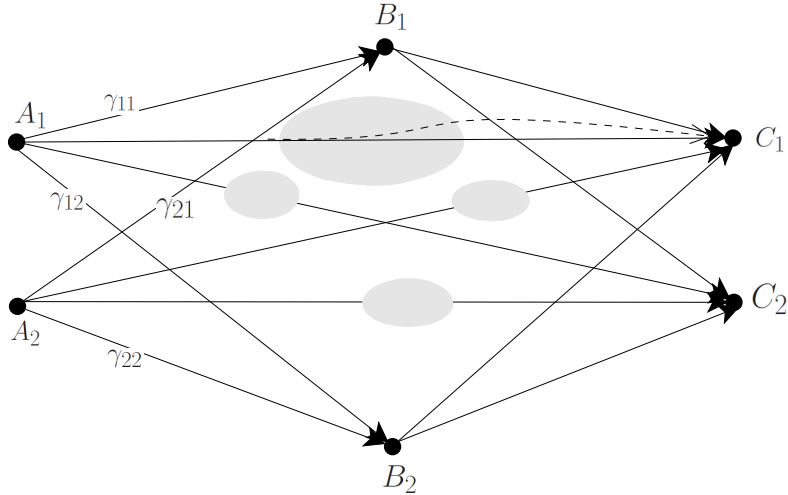}}
\caption{Points $A_i$, $B_i$, $C_i$, and the neighborhoods  $U^{ij}$ (gray ellipses) where we perturb the metric such that the  ``new'' geodesics (punctured line in the bigger gray ellipse) come  to the point $C_i$ from new directions. }\label{fig}
 \end{figure}

Let us construct, for each $j=1,...,n$, a mapping
$$\phi_{j}^{AB}: S^1_{A_1}  \times ...\times S^1_{A_n} \to S^1_{B_j}.$$

Denote by  $\gamma_{ij}=\gamma_{ij}^{AB}$ the geodesics connecting the $A$-points to the $B$-points,
we assume $\gamma_{ij}(0)= A_i$ and $\gamma_{ij}(1)=B_j$, see Fig.\ \ref{fig}. Now, for  an arbitrary  element
$(\alpha_1,...,\alpha_n)\in S^1_{A_i}  \times ...\times S^1_{A_n}$, consider the 1-form  $\beta_j\in S^1_{B_j}$
such that
 \begin{equation}  \label{condition}
\beta_j(\dot\gamma_{ij}(1))= \alpha_i(\dot\gamma_{ij}(0)), \ \ i= 1,...,n.
 \end{equation}
Clearly,  \eqref{condition} is a system of $n$ linear equations on $n$ components of the 1-form $\beta_j$.
The rows of the coefficient matrix of this system are the vectors $\dot\gamma_{ij}(1)$, $i=1,...,n$,
and the right hand side depends on $\dot\gamma_{ij}(0)$ and also on the 1-forms  $\alpha_1,..., \alpha_n$. Since
the vectors $\dot\gamma_{ij}(1)$, $i=1,...,n$, are linearly independent by our assumptions, the coefficient matrix is
nondegenerate and the system has precisely one solution $\beta_j$. We set $\phi^{AB}_j(\alpha_1,...,\alpha_n)= \beta_j$.

Define $\Phi^{AB}:S^1_{A_1}  \times ...\times S^1_{A_n} \to S^1_{B_1}  \times ...\times S^1_{B_n}$  by  setting  $$\Phi^{AB}(\alpha_1,...,\alpha_n) = \left(\phi^{AB}_1(\alpha_1,...,\alpha_n),...,\phi^{AB}_n(\alpha_1,...,\alpha_n)\right).$$

Similarly, we define $\Phi^{BC}:S^1_{B_1}  \times ...\times S^1_{B_n} \to S^1_{C_1}  \times ...\times S^1_{C_n}$,  $\Phi^{AC}$, $\Phi^{CA}$,  $\Phi^{CB}$. Clearly, these maps are  linear and actually they are isomorphisms since by construction $\Phi^{BA}$ is inverse to $\Phi^{AB}$ and so on.
Now consider
 \begin{equation} \label{eq}
\Phi^{BC}\circ \Phi^{AB}(\alpha_1,...,\alpha_n) = \Phi^{AC}(\alpha_1,...,\alpha_n).
 \end{equation}
We view this equality as a system of $n^2$ linear equations on $n^2$ components $(\alpha_1,...,\alpha_n)$ (recall that each $\alpha_i$ is a $1$-form and has therefore $n$ components).

Suppose a function $F$ which is linear in velocities, $F(x,v)= \sum_{i=1}^n K_i v^i$,  is an integral of the geodesic
flow. Denote by $\alpha_i\in S_{A_i}^1 $ the restriction of $F$ to $A_i$,
by $\beta_i\in S_{B_i}^1 $ the restriction of $F$ to $B_i$, and by $\delta_i\in S_{C_i}^1 $ the restriction of $F$ 
to $C_i$. Then equations \eqref{condition} defining the mappings $\Phi$ and also equation \eqref{eq} are fulfilled 
in view of \eqref{condkill}.

Thus, if there exists a nontrivial Killing 1-form, the system \eqref{eq} has a nontrivial solution. Our goal is therefore
to show that after a small $C^2$ perturbation the system (and therefore all systems in $\varepsilon'$-neighborhood
of it) has no nontrivial solutions.

The perturbation will be made in $n^2$ small nonintersecting neighborhoods $U^{ij}$, see Fig.~\ref{fig}. The neighborhood $U^{ij}$ intersects with no geodesics used to construct the mapping $\Phi^{AB}$, $\Phi^{BC}$, $\Phi^{AC}$ with exception of  the geodesic $\gamma^{AC}_{ij}$ that it does intersect. The perturbation is such that the initial velocity vector
$\dot \gamma^{AC}_{ij}(0)$ at the point $A_i$ remains unchanged, while the final velocity vector
$\dot \gamma^{AC}_{ij}(1)$ at the point $C_j$ is changed.
Clearly, such perturbation can be made $\varepsilon$-small in the $C^2$-topology, and clearly, by such a
perturbation, we can make the new vector $\dot \gamma^{AC}_{ij}(1)$ to be an arbitrary vector such that $g(\dot \gamma^{AC}_{ij}(1), \dot \gamma^{AC}_{ij}(1))= g(\dot \gamma^{AC}_{ij}(0), \dot \gamma^{AC}_{ij}(0))$
lying in a sufficiently small neighborhood of the old vector $\dot \gamma^{AC}_{ij}(1)$.

It is almost evident  that  by   such a small perturbation   one can achieve that the system \eqref{eq} has no nontrivial
solutions. Indeed, our perturbation does not change the mappings $\Phi^{AB}$ and $\Phi^{BC}$. But it does change the map
$\Phi^{AC}$ since $j$-rows of the matrix corresponding to $\phi_j^{AC}$ are the velocity vectors $\dot \gamma_{ij}(1)$, and the perturbation is small, but otherwise  is almost  arbitrary -- the only restriction is the $g$-length of
the rows of the matrix.
Thus  we can   force, by a small perturbation of the metric,  the system \eqref{eq} to have no nontrivial solution.
This perturbed metric, which we denote by $\tilde g$, does not have a nontrivial Killing 1-form.

Since a small perturbation of an unsolvable system of linear equations is also unsolvable, and a small $C^2$-perturbation of $\tilde g$ results in a small change of the geodesics, then for some $\varepsilon'>0$ there exists no metric
that is $\varepsilon'$-close to $\tilde g$ and has a nontrivial Killing 1-form.

Thus Theorem \ref{thm:1} is proven for $d=1$. \qed

\subsection{Proof for arbitrary $d$.} \label{sec:2.2}

We again  work locally, in a small geodesically convex disc $D\subset \mathbb{R}^n$ assuming $n\ge 2$.
Set $N=\binom{n+d-1}{d}$. It is known that $N$ is the dimension of the space $S^d$ of homogeneous polynomials
of degree $d$ in $n$ variables.

Take  $\kappa\in \mathbb{N}$ (which will be specified later; we'll see that it  is enough to take $\kappa=3$ but to keep the proof easier it is convenient to allow $\kappa$ to be big enough). Let us
consider $(\kappa+2)N$ points on the disc $D$ which we denote by
$A_1,...,A_N, B_{1,1}...,B_{1,N},...,  B_{\kappa,1},..., B_{\kappa,N},  C_1,...,C_{N}$.

As in the case $d=1$, we assume that no three of the points $ A_1,...,C_{N}$ lie on one geodesic.
In order to introduce the second assumption on the points, we call a set of $N$
vectors of  $\mathbb{R}^n$\ $d$-\emph{decisive}, if the values of any homogeneous polynomial of degree $d$  on this set
determine  this polynomial. For example,
for $d=1$ we have $N=n$ and $1$-decisive sets are precisely those containing $n$ linearly independent vectors.
Clearly,  the  set of $d$-decisive sets is open and everywhere dense in
the set of all subsets of $\mathbb{R}^n$ containing $N$ elements.

As in the case $d=1$, we denote the set $\{A_1,...,A_N\}$ by $A$, the  set $\{B_{\ell,1},...,B_{\ell,N}\}$ by $B_\ell$,
and the set $\{C_1,...,C_{N}\}$ by $C$.
Assume that for every $M_1\neq M_2\in \{A,B_1,...,B_\kappa, C\}$ the initial vectors of geodesics connecting
any point of $M_1$  to all points of $M_2$ form a $d$-decisive set in the tangent space to this point of $M_1$.

Denote by $S^d_{A_i}$ (resp.\ $S^d_{B_{\ell,i}}$, $S^d_{C_i}$) the space of  homogeneous polynomials of degree $d$
on $T_{A_i}D$ (resp.\ on $T_{B_{\ell,i}}D$, $T_{C_i}D$).
These are $N$-dimensional vector spaces. We will denote the elements of
$S^d_{A_i}$ by $\alpha_i$ and the elements of $S^d_{B_{\ell,i}}$ by $\beta_{\ell,i}$.

Let us construct, for each $j=1,...,N$, and $\ell=1,...,\kappa$, a mapping
 $$
\phi_{\ell,j}^{AB}: S^d_{A_1}  \times ...\times S^d_{A_N} \to S^d_{B_{\ell,j}}.
 $$

Denote by $\gamma^\ell_{ij}$ the geodesics connecting $A$-points to $B$-points:
$\gamma^\ell_{ij}(0)= A_i$, $\gamma^\ell_{ij}(1)=B_{\ell,j}$. For an arbitrary element
$(\alpha_1,...,\alpha_N)\in S^d_{A_1} \times ...\times S^d_{A_N}$
find an element $\beta_{\ell,j}\in S^d_{B_{\ell,j}}$ such that
 \begin{equation}  \label{conditionN}
\beta_{\ell,j}(\dot\gamma^\ell_{ij}(1))= \alpha_i(\dot\gamma^\ell_{ij}(0)), \ \ i= 1,...,N.
 \end{equation}
Clearly,  \eqref{conditionN} is a system of $N$ linear equations on the $N$ components of $\beta_{\ell,j}$.
Since the vectors $\dot\gamma_{ij}(1)$, $i=1,...,N$, form a $d$-decisive set in $T_{B_j}D$, the system has
precisely one solution $\beta_{\ell,j}$. We set $\phi^{AB}_{\ell,j}(\alpha_1,...,\alpha_N)= \beta_{\ell,j}$.

Define $\Phi_{\ell}^{AB}:S^d_{A_1}  \times ...\times S^d_{A_N} \to S^d_{B_{\ell,1}}  \times ...\times S^d_{B_{\ell,N}}$  by  setting  \
 $$
\Phi_{\ell}^{AB}(\alpha_1,...,\alpha_N) = \left(\phi^{AB}_{\ell,1}(\alpha_1,...,\alpha_N),...,\phi^{AB}_{\ell, N}(\alpha_1,...,\alpha_N)\right).
 $$

Similarly, we define $\Phi_\ell^{BC}:S^d_{B_{\ell, 1}} \times ...\times S^d_{B_{\ell,N}} \to S^d_{C_1} \times ...\times S^d_{C_N}$. Clearly, these maps are linear isomorphisms.

Now consider, for any $\ell \in\{2,...,\kappa\}$, the relation
 \begin{equation} \label{eqN}
\Phi_\ell^{BC}\circ \Phi_\ell^{AB}(\alpha_1,...,\alpha_N) = \Phi_1^{BC}\circ \Phi_1^{AB}(\alpha_1,...,\alpha_N).
 \end{equation}
We view this equality as  a system of $N^2$ linear equations on $N^2$ components $(\alpha_1,...,\alpha_N)$ (recall that each $\alpha_i$ is a homogeneous polynomial of degree $d$ and has $N$ components). The coefficients in this system of equations are explicit algebraic functions  in the velocity vectors of the geodesic segments
connecting $A$-points to $B$-points and $B$-points to $C$-points.

Suppose a function $F$, which is a homogeneous polynomial in velocities of degree $d$,
is an integral of the geodesic flow of $g$.
Denote by $\alpha_i\in S_{A_i}^d $ the restriction of $F$ to $A_i$, and by $\beta_{\ell, i}\in S_{B_{\ell,i}}^d $
the restriction of $F$ to $B_{\ell,i}$. Then the equations \eqref{conditionN} defining the mappings $\Phi^{AB}_\ell$
are fulfilled together with the corresponding equation for $\Phi^{BC}_\ell$ in view of \eqref{condkill}.
Consequently equations \eqref{eqN} are also fulfilled.

 \begin{remark}
We have chosen a different approach with \eqref{eqN} for $d>1$ than that with \eqref{eq} for $d=1$.
It is also possible to proceed with the latter for $d>1$ by examining the system of equations $\Phi_\ell^{BC}\circ\Phi_\ell^{AB}=\Phi^{AC}$ for various $\ell$ but the arguments become more cumbersome.
 \end{remark}

We conclude  that if system \eqref{eqN} has only trivial solutions, the metric admits no nontrivial integral
that is polynomial in velocities of degree $d$.
Notice that the notion of the triviality of a solution depends on the parity of the number $d$. For odd $d$ the trivial
solution is $\alpha_1=...=\alpha_N=0$: it corresponds to the identically zero integral. For even $d=2q$ we treat
the solutions $\alpha_1=c\cdot H^q_{|A_1},...,\alpha_N= c\cdot H^q_{|A_N}$ ($c\in\mathbb{R}$) as trivial:
such solutions correspond to the $q$-th power of the Hamiltonian $H(x,v)=\tfrac12\sum_{i,j}g_{ij}(x)v^iv^j$.

\begin{figure}[ht!]
{\includegraphics[width=0.7\textwidth]{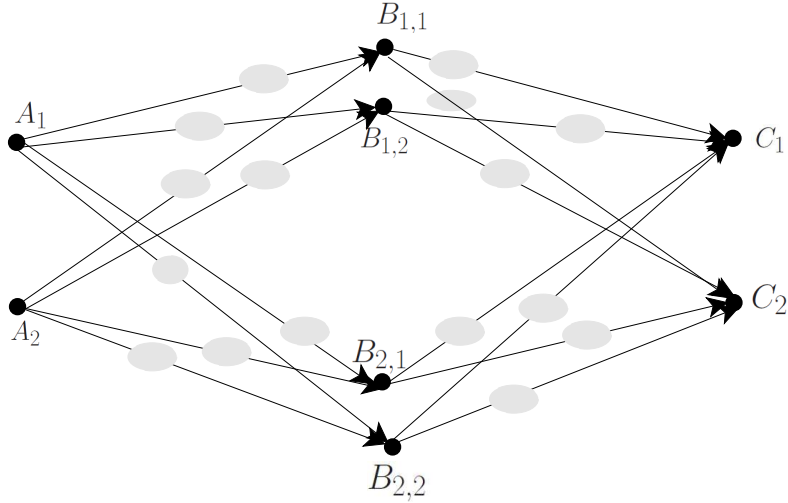}}
\caption{The points $A_i$, $B_{\ell,i}$, $C_i$, the geodesics used in equation \eqref{eqN}, and  the  neighborhoods
(gray ellipses) supporting perturbations of the metric to achieve different velocity vectors at endpoints of the
new geodesics.}\label{fig1}
 \end{figure}

We shall now show that for a given  metric $g$ there exists  a small $C^2$-perturbation $\tilde g$ of it  such that
\eqref{eqN} has only trivial solutions. This perturbation will be localized in small neighborhoods
of points of geodesics connecting $A$-points to $B$-points and $B$-points to $C$-points. As in the case $d=1$,
these neighborhoods  will be chosen in such a way that each neighborhood intersects with only  one of the geodesics used to construct the mappings $\Phi_\ell^{AB}$, $\Phi_\ell^{BC}$, see Fig.\ \ref{fig1} and the  caption to it.
With such a perturbation we can achieve any (sufficiently small)  change of the velocity vectors at the initial
and the final points of these geodesics, with the only restriction that for any geodesic
the length of the velocity vectors at the initial and the final points remains the same.

If after such a small perturbation system \eqref{eqN} has only the trivial solutions then we are done, because
a small $C^2$-perturbation of the metric implies a small perturbation of the coefficients in system \eqref{eqN}
and cannot increase the dimension of the solution space of \eqref{eqN}, which is clearly an upper semi-continuous function. In particular, if for $\tilde g$ the system has only trivial solutions,
so will it be for its small  $C^2$-perturbations.

Suppose, by the way of proof ad absurdum, that for any sufficiently small perturbation the system \eqref{eqN} has
nontrivial solutions. Because the coefficient matrix of system \eqref{eqN} depends algebraically on the velocity
vectors at the initial and final points of the geodesic segments used in the construction, this system will
have a nontrivial solution for any initial and final vectors used instead of $\dot\gamma^\ell_{ij}(0)$ and $\dot\gamma^\ell_{ij}(1)$ in formula \eqref{conditionN}.

In what follows we are going to show that, by replacing the initial and final vectors of geodesics by
appropriately chosen vectors from a decisive set, one can achieve that the mappings
$\Phi_{\ell}^{BC}\circ\Phi_\ell^{AB}$ correspond to some orthogonal linear transformations $\sigma_\ell\in O(g)$
(the notion $O(g)$ is well-defined since we work in coordinates such that at all $A$-, $B$-, $C$-points
the metric is given by the same matrix, and our perturbation does not change $g$ at the $A$-, $B$-, $C$-points).
Let us emphasize once more that we have the following two restrictions on the \emph{replacements}:
for any $j$ the replacement for the set $\dot \gamma^\ell_{ij}(1)$ (used in the construction of the mapping
$\phi^{AB}_{\ell} \circ\phi^{BC}_{\ell}$) should form a decisive set, and  for any such a geodesic
$\gamma^\ell_{ij}$ we replace its initial and final velocity vectors $\dot \gamma^\ell_{ij}(0)$ and
$\dot \gamma^\ell_{ij}(1)$ by vectors of the same length. Besides these restriction, the replacement can be arbitrary
and replacements with this restriction will be called \emph{admissible}.

At the tangent space to each of the points $A_1,...,C_N$ choose a coordinate system such that the matrix of $g$ is the
standard (e.g.\ the identity matrix if $g$ has Riemannian signature and $\operatorname{diag}(-1,1,...,1)$ if $g$ has
Lorenzian signature). This gives us a coordinate system on the spaces $S^d_{A_i}$, $S^d_{B_{\ell,i}}$, $S^d_{C_i}$.
In this coordinate system, for a fixed $\ell\in \{1,...,\kappa\}$, the linear  mapping $\Phi_\ell^{BC}\circ\Phi_\ell^{AB}$ can be viewed as an $N^2\times N^2$-matrix. Let us show that we can admissibly replace the vectors
$\dot \gamma^\ell_{ij}(0)$,  $\dot \gamma^\ell_{ij}(1)$ such that this matrix has a relatively simple form,
in particular it is block-diagonal with $N$ blocks of dimension $N\times N$ and all these blocks are the same.

Let us first explain that one can make, for $\ell=1$,  all these $N\times N$-blocks to be equal the identity matrix
(which of course makes the entire $N^2\times N^2$-matrix the identity matrix).
Choose and fix a decisive set $v_1,...,v_N\in \mathbb{R}^n$ such that $g(v_i, v_i)= 1$ at $A_1$  and replace, 
in the construction of $\Phi_1^{AB}$, the vectors $\dot \gamma^1_{ij}(0)$ by $v_i$ and the vectors 
$\dot \gamma^1 _{ij}(1)$ by $v_j$. The replacement is admissible since by our choice of the coordinates at  
the $A$-, $B$- and $C$-points we have  $g(v_i, v_i)= 1$ at all these points.
After this replacement  formula \eqref{conditionN} becomes $\beta_{1,j}(v_i)=\alpha_i(v_j)$.
Similarly, replace in the construction of $\Phi_1^{BC}$, the vectors $\dot\gamma^1_{ij}(0)$ by $v_j$ and  the vectors
$\dot \gamma^1_{ij}(1)$ by $v_i$. After this replacement, the   formula defining  $\Phi_1^{BC}$ reads
$\delta_j(v_i)=\beta_{1,i}(v_j)$, where $\delta_i$ denote the elements of $S^d_{C_i}$. We see that
the equations used in the construction of $\Phi_1^{AB}$ are dual to that used in the construction of $\Phi_1^{BC}$
and therefore the mapping $\Phi_1^{BC}\circ \Phi_1^{AB}$ is given by the identity $N^2\times N^2$ matrix, as we wanted.

Let us consider  $\ell\ge 2$. For any element $\sigma_\ell\in O(g)$, let us construct an admissible  replacement such that each $N\times N$-block  of the matrix of $\Phi_\ell^{BC}\circ \Phi_\ell^{AB}$ corresponds to the matrix of the
linear transformation induced on $S^d_{A_i}$ by $\sigma_\ell$, viewed as a linear transformation of $T_{A_i}D$.
This replacement does not change the mappings  $\Phi_m^{BC}\circ \Phi_m^{AB}$ for $m\ne \ell$, because we will
replace only the vectors $\dot\gamma^\ell_{ij}(0)$, $\dot\gamma^\ell_{ij}(1)$ and these vectors are not used 
in the construction of the mappings $\Phi_m^{BC}\circ \Phi_m^{AB}$.

To do so, we  replace, in the construction of $\Phi_\ell^{AB}$, the vectors $\dot \gamma^\ell_{ij}(0)$ by $v_i$
and the vectors $\dot \gamma^\ell_{ij}(1)$ by $v_j$. Equations \eqref{conditionN} become
$\beta_{\ell,j}(v_i)=\alpha_i(v_j)$. Now, slightly different from above,
in the construction of $\Phi_\ell^{BC}$, we replace the vectors $\dot\gamma^\ell_{ij}(0)$ by $v_j$ and the vectors
$\dot \gamma^\ell_{ij}(1)$ by $\sigma_\ell(v_i)$. Then, equations \eqref{conditionN} become
$\delta_j(\sigma_\ell(v_i))=\beta_{\ell,i}(v_j)$. Clearly, after this change the matrix of $\Phi_\ell^{BC}\circ \Phi_\ell^{AB}$ is $\operatorname{diag}(S^d\sigma_\ell,...,S^d\sigma_\ell)$, as we want.

By construction, any solution $\alpha_1,...,\alpha_N$ of system \eqref{eqN} should be invariant
with respect to the action induced by $\sigma_\ell$, i.e., each homogeneous polynomial $\alpha_i$ should be invariant:
$\alpha_i= \sigma_\ell^*\alpha_i$.

Now we specify $\kappa$: choose elements $\sigma_2,...,\sigma_{\kappa}\in O(g)$ such that any homogeneous
polynomial of degree $d$ in $\mathbb{R}^n$, which is invariant under all $\sigma_\ell$, is trivial.
Let us explain why we can choose a finite set of such elements. The condition $\sigma^*\alpha= \alpha$
is a system of linear equations on $\alpha$ whose coefficients come from $\sigma\in O(g)$.
Consequently, the condition that a homogeneous polynomial $\alpha$ of degree $d$ is invariant with respect
to the whole group $O(g)$ is the following infinite system of linear equations:
 $$
\sigma^*\alpha=\alpha,  \ \textrm{where $\sigma  \in O(g)$}.
 $$
Since any homogeneous polynomial of degree $d$ invariant with respect to the whole group $O(g)$ is trivial,
we obtain that the rank of this system of equations is $N$ for $d$ odd and $N-1$ for $d$ even. This implies that
one can choose a finite subsystem of this system such that its rank is still $N$ for $d$ odd and $N-1$ for $d$ even.

Thus, with this $\kappa$ and the choice $\sigma_2,..., \sigma_{\kappa}\in O(g)$ we obtain that system \eqref{eqN}
has only trivial solutions, as we want. Therefore Theorem \ref{thm:1} is proven. \qed

 \begin{remark}
We can minimize the value of $\kappa$ in the above proof.
Indeed, for odd $d$ we simply take $\kappa=2$ and $\sigma_2= -\operatorname{id}$, since any homogeneous polynomial
of odd degree, which is invariant under the reflection with respect to the zero point, is identically zero.
For even $d$, $\kappa=3$ is enough because $O(g)$ has a dense free subgroup generated by two elements
$\sigma_2,\sigma_3$, as follows from \cite[Theorem 1.1]{BG}.
 \end{remark}

\subsection{Proof of Corollary \ref{cor:0}.} \label{sec:2.3}

In the proof of Theorem \ref{thm:1} we have shown  how to reconstruct the coefficients of a polynomial in velocities
integral of degree $d$ by its coefficients in $N= \binom{n+d-1}{d}$ points in general position. The reconstruction
requires only the exponential mapping from these points and is algebraic otherwise. Since for a $C^k$-smooth metric the
exponential mapping (viewed as a local flow on $TD$) is $C^{k-1}$-smooth by standard (see e.g. \cite[\S V.3]{hartman}) results on smooth dependence of solutions of ODE on initial conditions and parameters, the reconstruction gives us a $C^{k-1}$-smooth integral. \qed

\section{Proof of Theorem \ref{thm:2}.} \label{sec:2.4}

We will first recall the following statement which is well known in the folklore.

 \begin{theorem} \label{folklore}
Let $g$ be a metric of smoothness $C^{d+2}$ on $M$. Then there exists a 
vector bundle of rank $\tfrac{(n+d-1)!\,(n+d)!}{(n-1)!\,n!\,d!\,(d+1)!}$ over $M$ with a linear connection
canonically constructed by $g$ such that parallel sections of this connection are in one-to-one correspondence with
integrals for the geodesic flow of $g$ that are polynomial in momenta of degree $d$.

In local coordinates, the coefficients of the connection (Christoffel symbols) are   rational functions of
the components of the metric $g$ and in their partial derivatives up to order $d+1$.
 \end{theorem}

Theorem \ref{folklore} mostly follows from \cite{Thompson} and was proved but not stated in \cite[\S 4]{Wolf}.
Let us briefly explain it (because in \cite{Thompson}  the constant curvature case was studied, even though it is not
crucial). The original Killing equation is \eqref{killing}, and its $k$-th prolongations are given by
 \begin{equation}\label{ProlKill}
K_{(i_1...i_d, i_{d+1}\!)j_1...j_k}\equiv 0.
 \end{equation}
Since the symbol of $K_{i_1...i_d, i_{d+1}j_1...j_k}$ is symmetric by the last $(k+1)$ indices,
the symbol of the equation describes the kernel of the map $S^d\otimes S^{k+1}\to S^{d+1}\otimes S^k$
(symmetrization by the first $(d+1)$ indices).
It is easy to check that this map is epimorphic for $k\leq d$, which means that there are no relations between
the linear equations of the system (compatibility conditions) up to this prolongation. Moreover for $k=d$
this map is an isomorphism, which means that the Cartan distribution $\Delta$ of the equation
$\mathcal{E}_{d+1}\subset J^{d+1}$, considered as a subset in the space of $(d+1)$-jets given by \eqref{ProlKill}
for all $k\leq d$ and all possible indices, is transversal to the fibers and the natural projection
$\pi:J^{d+1}\to M$ maps
$d\pi:\Delta_{x_{d+1}}\stackrel{\sim\,}\to T_xM$ for every $x_{d+1}\in\mathcal{E}_{d+1}\cap\pi^{-1}(x)$.
This $(\mathcal{E}_{d+1},\pi,M)$ is the required vector bundle, and
$\Delta$ is the required linear connection on it. Solutions to the Killing equations are bijective (via the
jet-prolongation) to the integral surfaces of $\Delta$.

Now, by Theorem \ref{thm:1}, for any dimension  $n\ge 2$  and for any $d$,  there exists a real analytic metric $\hat g$  on a disc $D\subseteq \mathbb{R}^n$ such that the geodesic flow of the   metric does not admits a nontrivial  integral that is polynomial in momenta of degree $d$. For this metric, the connection on the bundle from Theorem \ref{folklore}, which we denote by $\nabla$,
is also real analytic.

We consider the curvature of this connection. In a local coordinate system, it has components
$R^\alpha_{\ \beta ij}$, where $i,j\in \{1,...,n\}$ and $\alpha, \beta\in \left\{1,..., \tfrac{(n+d-1)!(n+d)!}{(n-1)!n!d!(d+1)!}\right\}.$

It is known (see \cite[Theorem 5]{ozeki}, in fact it follows immediately from the Ambrose-Singer Theorem) that
for a linear real analytic connection the holonomy algebra of the connection at point $x_0$ is generated by 
the covariant derivatives, up to a certain order, of the curvature of the connection. 
The parallel sections of a connection correspond to elements of the fiber at $x_0$ that
are invariant with respect to the holonomy group (and so annihilated by the holonomy algebra). 
Since there exists no nontrivial parallel section,
there exists $N \in \mathbb{N}$ such that no nontrivial element $V=(V_\alpha)$
of the fiber at $x_0$ satisfies, for all $k\in \{0,...,N\}$ and for any  
$i,j, i_1,...,i_k\in \{1,...,n\}$, the  linear equations (where $\nabla_i=\nabla_{\partial_{x^i}}$):
 \begin{equation} \label{last}
\sum_\alpha V_\alpha\cdot\left(\nabla_{i_1}\cdots\nabla_{i_k}\, R^\alpha_{\ \beta ij}\right)=0.
 \end{equation}

Let us now consider, for an arbitrary metric $g$ on $D$, the following 1-parameter family
 $$
g_t= t \hat g  + (1-t)  g.
 $$
In the Riemannian case all $g_t$ are metrics; however in the general case we can only ensure that
$g_t$ are metrics for $t$ close to $1$ and for $t$ close to $0$,
and the elements of this family for other $t$ are not important for us (for any point $x\in M$ at most $n$ 
values of $t\in[0,1]$ correspond to degenerate 2-tensor $g_t$ at $x$, so for a generic $t$ the tensor $g_t$
is a metric in a small neighborhood of $x$).

Consider the system of linear  equations of type \eqref{last} but constructed by the  metric $g_t$. The system   is  fulfilled for any parallel section, so if the system of equations has only trivial solution then any parallel section is trivial. The coefficients of this system depend algebraically on $t$, and on $(N+d+2)$-jets of the metrics $g$ and $\hat g$
at $x_0$. Since system \eqref{last} has no nontrivial solution for $t=1$, it has no nontrivial solutions
for almost all $t$; thus, there exists arbitrary small $t>0$ such that the system \eqref{last} for $\tilde g= g_t$ has
no nontrivial solution. Moreover, for any metric $g'$ such that its  $(N+d+2)$-jet at the point $x_0$ is  close
to that of  $\tilde g$, system \eqref{last} for $g'$ is obtained by a small perturbation of that for $\tilde g$ and again has no nontrivial solutions, as we want.  Theorem \ref{thm:2} is proven. \qed

 \begin{remark}
In the proof of Theorem \ref{thm:2}, in order to apply Theorem  \ref{folklore}, we only need the existence of one
real-analytic metric (on a disc of an arbitrary dimension) that does not admit a nontivial  integral of degree $d$. Theorem \ref{thm:1} gives as such a metric in  any dimension.   Let us note that in dimension 2 the  existence of such a  metric has been known before, see e.g.\ \cite{Ten}. In fact, in order to apply
Theorem  \ref{folklore} in dimension $n$ it suffices to have the existence of one real analytic metric  on a simply connected $n$-manifold, such that the geodesic flow of this metric does not admit a nontrivial  integral of degree $d$,  and again in dimension 2 the existence of such a metric on the sphere  follows from results of \cite{Bolotin,knieper,paternain,taimanov}, as we noted in the Introduction.
 \end{remark}

\section{Conclusion: local smooth integrability of geodesic flows \\ and other classes of integrals. } \label{sec:2.5}

In our paper we have demonstrated generic non-existence of integrals analytic in momenta (and smooth otherwise).
Let us show now the local existence of integrals that are smooth in momenta. This folklore known result is the
reason the general smooth integrability problem we mentioned in the introduction has sense only globally.

 \begin{proposition} \label{lastt}
Let $D\subset \mathbb{R}^n  $ be a geodesically convex  disc  for a metric  $g$ on $\mathbb{R}^n$ of arbitary signature. Then the geodesic flow of $g|_D$ is Liouville-integrable.
 \end{proposition}

By \emph{Liouville  integrability} we understand the existence of $n$ functionally independent almost everywhere integrals in the involution. From the proof it will be clear that actually it is also super-integrable, i.e., there exists $(2n-1)$ functionally independent integrals such that $n$ of them are in the involution.

We will see in the proof that actually the only condition we need is that the symplectic reduction
of the hypersurface in $T^*D$ given by $|H|=1$ is a smooth manifold. 
In the case when $D$ is geodesically convex this property is evident,
but it also holds in the case of, for example, simple Riemannian manifolds in the terminology of \cite{preprint} 
and follows from  \cite[Theorem 1.3]{preprint} for dimension $n=2$.

 {\bf Proof.}
We discuss first the Riemannian case, where the arguments are slightly simpler. We consider $T^*D$ with the canonical symplectic structure and the Hamiltonian action on it generated by the Hamiltonian $H$ of the geodesic flow. Of course, the action is only locally defined and the projection of its orbits to $D$ are geodesic segments; because of geodesic convexity the endpoints of these segments lie on the boundary $\partial D^n\simeq S^{n-1}$.

Take the energy level $E_1= \{(x,p)\in T^*DF \mid H(x,p)=1\}$ and consider the symplectic reduction corresponding
to this level. The result is a symplectic manifold of dimension $2(n-1)$, which we denote by $Q$. Though it is not important for our proof, let us mention that in the Riemannian case $Q$ is diffeomorphic to $TS^{n-1}$. Indeed,
the points of $Q$ are oriented geodesic segments with endpoints on the boundary, and so can be parameterised by
outward unit vectors along the boundary (whose natural projection to $TS^{n-1}$ gives an open neighborhood
of the zero section).

Denote by $\pi:E_1\to Q$ the tautological projection. It is a general property of the symplectic reduction that
for any function $F\in C^\infty(Q)$ the function $\pi^*F\in C^\infty(E_1)$ is an integral of the restriction of the
geodesic flow to $E_1$.

Now take $(n-1)$ smooth functionally independent almost everywhere Poisson-commuting functions
$F_2,...,F_n$ on $Q$. The existence of such is standard, see e.g.\ \cite[Proposition 1 in \S5.1]{fomenko}.
In addition to involutivity, following the proof of \cite[Proposition 1 in \S5.1]{fomenko}  we can arrange that these functions and all their partial derivatives are bounded
(each by its own constant). We extend them to $T^*D$ by the formula
$\hat F_i=\exp(-H^{-2})\cdot\nu^*(\pi^*F_i)$, where $\nu:T^*D\setminus \{0\} \to E_1$
is the radial projection $(x,p)\mapsto\left(x,\frac{p}{\sqrt{H}}\right)$.  We set $\hat F_i(x,0)=0$.
It is straightforward to check that the functions $\hat F_i$ are smooth involutive functionally independent
integrals for $H$, which proves Proposition \ref{lastt} for Riemannian metrics. 

Consider now the metrics of indefinite signature. The energy levels $E_{\pm}=\{(x,p)\in T^*DF: H(x,p)=\pm1\}$ admit symplectic reductions to symplectic manifolds $Q_{\pm}$. Take $(n-1)$ smooth functionally independent almost everywhere Poisson-commuting functions $F_2,...,F_n$ on $Q_+\cup Q_-$, which are bounded as well as all their derivatives.
Extend them to $T^*D$ as before: for $(x,p)$ with $H(x,p) \neq 0$ set
$\hat F_i=\exp(-H^{-2})\cdot\nu^* \left(\pi^* F_i  \right)$,
where $\pi:E_\pm\to Q_\pm$ is the tautological projection and $\nu: (x,p)\mapsto\left(x,\frac{p}{\sqrt{|H|}}\right)$ is the radial projection. For $(x,p)$ with $H(x,p) =0$ set $\hat F_i(x,p)=0$. We obtain  smooth involutive functionally independent integrals $\hat F_1=H,\hat F_2,...,\hat F_n$ on $T^*D$, and Proposition \ref{lastt} is proved. \qed

Note that even if the metric is real-analytic (which makes the manifolds $Q$ and $Q_\pm$
from the proof of Proposition real-analytic and after investing some work one also obtains real-analytic integrals
$\tilde F_i$ on $E_1$ or on $E_{\pm}$, our proof does not give real-analytic integrals $\hat F_i$, since the
conformal factor $\exp(H^{-2})$ (extended by $0$ to $\{H=0\}$) is not real analytic around the zero section.
Actually our results, in particular Corollary \ref{cor:2}, explain that for generic metric $g$ the integrals do
not real-analytically extend from $E_1$ (or from $E_\pm$) to the entire $T^*D$.

Let us also discuss some other classes of integrals. Examining the proof of Theorem \ref{thm:1} we see that
the only assumption we used is that the restriction of an integral to each tangent space is determined by a
finitely many numbers (coefficients $K_{i_1...i_d}$ in the case of integrals polynomial in momenta).
Therefore the results remains true if we replace the integrals polynomial in momenta
by any other class, satisfying this property.
In particular, this holds for integrals that are rational in momenta:
 $$
F(x,p)=\frac{\sum_{i_1,...,i_d=1}^n K_{i_1...i_d}(x)\,p_{i_1}...p_{i_d}}{\sum_{i_1,...,i_r=1}^n Q_{i_1...i_r}(x)\, p_{i_1}...p_{i_r}}.
 $$
Such integrals were studied already in \cite{Darboux} and the question of their existence and nonexistence,
for real-analytic metrics on surfaces, is the subject of the recent paper \cite{kozlov}.


\begin{thebibliography}{100}
  	
\bibitem{Bolotin}
S. Bolotin, P. Negrini, \emph{A variational criterion for nonintegrability}, Russian J. Math. Phys. {\bf 5} (1997), 415-436.

\bibitem{BG}
E. Breuillard, T. Gelander, {\it On dense free subgroups of Lie groups}, Journal of Algebra {\bf 261} (2002), no.2, 448-467.

\bibitem{burns} K. Burns, V. S. Matveev, \emph{ Open problems and questions about geodesics,} arXiv:1308.5417.

\bibitem{Darboux}
G. Darboux, \emph{Le\c{c}ons sur la th\'eorie g\'en\'erale des surfaces},
Vol. III, Chelsea Publishing, 1896.

\bibitem{duna}
R.\,L. Bryant, M. Dunajski, M. Eastwood, \emph{Metrisability of two-dimensional projective structures},
Journ. Differential Geom. {\bf 83} (2009), no.3, 465-500.

\bibitem{Contreras}
G. Contreras, \emph{Geodesic flows with positive topological entropy, twist maps and hyperbolicity},
Ann. of Math. (2) {\bf 172} (2010), no. 2, 761-808.

\bibitem{ContrerasPa}
G. Contreras, G. Paternain, {\it Genericity of geodesic flows with positive topological entropy on $S^2$},
Journ. Differential Geom. {\bf 61} (2002), no. 1, 1-49.

\bibitem{fomenko} A. T. Fomenko,  {\it Symplectic geometry}, Advanced Studies in Contemporary Mathematics 
{\bf 5} (1988), Gordon and Breach Science Publishers, New York.

\bibitem{hartman} P. Hartman, \emph{Ordinary differential equations}, John Wiley \& Sons, Inc., 
New York-London-Sydney, 1964.

\bibitem{knieper}
G. Knieper, H. Weiss, {\it $C^\infty$ genericity of positive topological entropy for geodesic flows on $S^2$},
Journ. Differential Geom. {\bf 62} (2002), no. 1, 127-141.

\bibitem{kozlov} V. V. Kozlov,  {\it On rational integrals of geodesic flows},
Regul. Chaotic Dyn. {\bf 19} (2014), no. 6, 601--606.

\bibitem{kruglikov}
B. Kruglikov, \emph{Invariant characterization of Liouville metrics and polynomial integrals,}
J. Geom. Phys. {\bf 58} (2008), no. 8, 979-995.

\bibitem{mathoverflow}
\textsf{http://mathoverflow.net/questions/178761/generic-absence-of-non-trivial-first-integrals-of-geodesic-flows}
(2014).

\bibitem{ozeki} H.  Ozeki, \emph{Infinitesimal holonomy groups of bundle connections,} Nagoya Math. J. {\bf 10} (1956), 105--123.
 
\bibitem{paternain}
G. Paternain, {\it Entropy and completely integrable Hamiltonian systems},
Proc. Amer. Math. Soc. {\bf 113} (1991), no. 3, 871-873.

\bibitem{preprint} G. Paternain, H. Zhou,  {\it Invariant distributions and tensor tomography on simple manifolds,} 
preprint.

\bibitem{taimanov}
I. Taimanov, {\it Topology of Riemannian manifolds with integrable geodesic flows}, Proc. Steklov Inst. Math.
{\bf 205} (1995), no. 4, 139-150.

\bibitem{Ten}
V.\,V. Ten, \emph{The Local Integrals of Geodesic Flows}, Regular Chaotic Dynam. {\bf 2} (1997),  87-89.

\bibitem{Thompson}
G. Thompson, \emph{Killing tensors in spaces of constant curvature}, J. Math. Phys. {\bf 27} (1986), no. 11, 2693-2699.

\bibitem{Wh}
E.\,T. Whittaker, {\it A Treatise on the Analytical Dynamics of Particles and Rigid Bodies}.
Cambridge University Press, Cambridge, 1937.

\bibitem{Wolf}
Th. Wolf, {\it Structural equations for Killing tensors of arbitrary rank}, Comput. Phys. Comm. {\bf 115} (1998),   316-329.

\end{thebibliography}
\end{document}